\documentclass[11pt]{article}
\def\diff{\mathsf{d}}
\usepackage[all]{xy}
\usepackage{color}
\usepackage{amscd}      
\usepackage{amssymb}
\usepackage[intlimits]{amsmath}

\usepackage{theorem}

\newtheorem{prop}{Proposition}[section]
\newtheorem{lem}[prop]{Lemma}

\newtheorem{coro}[prop]{Corollary}
\newtheorem{theo}[prop]{Theorem}

\theorembodyfont{\upshape}

\newtheorem{defi}[prop]{Definition}

\newtheorem{numrmk}[prop]{Remark}

\newtheorem{rmk}{Remark}

\newenvironment{pf}{\begin{trivlist}\item[]{\sc Proof.}}%
           {\nolinebreak $\Box$ \end{trivlist}}

\newcommand{\FF}{\mathcal F}
\newcommand{\superF}[1]{\FF_{#1}}

\newcommand{\R}{\mathbb R}

\newcommand{\N}{\mathbb N}

\newcommand{\ZZZ}{Z}

\newcommand{\id}{\rm id}
\newcommand{\ad}{\rm ad}
\newcommand{\lb}[1]{\left[ #1 \right]}

\newcommand{\bb}[1]{\left\{ #1 \right\}}
\newcommand{\LB}{\left[\cdot , \cdot \right]}

\newcommand{\partialt}{\partial^t}

\begin{document}

\title{{\bf Jacobi structures in supergeometric formalism.}}
\author{Paulo dos Santos Antunes, Camille Laurent-Gengoux,
\\[5pt] {\small\it Center for Mathematics, Universidade de Coimbra}
\\[5pt]
{\small\it e-mail: pantunes@mat.uc.pt,
claurent@mat.uc.pt} }
\date{}

\sloppy \maketitle

\begin{abstract}
We use the supergeometric formalism, more precisely, the so-called ``big bracket'' (for which brackets and anchors are encoded by functions on some graded symplectic manifold) to address the theory of Jacobi algebroids and bialgebroids (following mainly Iglesias-Marrero \cite{IM} and Grabowski-Marmo \cite{GM} as a guideline). This formalism is in particular efficient to define the Jacobi-Gerstenhaber algebra structure associated to a Jacobi algebroid, to define its Poissonization, and to express the compatibility condition defining Jacobi bialgebroids. Also, we claim that this supergeometric language gives a simple description of the Jacobi bialgebroid associated to Jacobi structures, and conversely, of the Jacobi structure associated to Jacobi bialgebroid.
\end{abstract}

\section{Introduction.}

  Assume that we give ourself a Poisson structure ${\mathcal P}$ on a manifold $X$ and an Euler vector field (i.e. a vector field ${\mathcal E}$ on $X$ such that ${\mathcal L}_{\mathcal E}{\mathcal P}={\mathcal P}$). A Lie-Poisson structure, taking for ${\mathcal E}$ the usual Euler vector field, is an example of that situation. More generally, any weight homogeneous Poisson structure \cite{LPV}, of weight $k \neq 2$, together with the weighted Euler vector field divided by $2-k$ gives an example. In this situation, for every hyper-surface $M \subset X$ transversal to ${\mathcal E}$, i.e such that
 $$T_m M \oplus {\mathbb R} {\mathcal E}_m = T_m X \quad \forall m \in M ,$$
there exists, on the submanifold $M$, an unique bivector field $\pi$ and an unique vector field $E$ such that:
  $$ {\mathcal P}_m = \pi_m +  {\mathcal E}_m \wedge E_m  \quad \forall m \in M.$$
 It turns that  the pair  $(\pi,E)$ satisfies the following two relations,
$$ \lb{\pi,\pi}= -2 E \pi \hbox{ and } \lb{E,\pi} =0.$$
 A bivector field and a vector field on a given manifold $M$ satisfying these conditions form what is called a Jacobi structure on $M$. It can be shown that any Jacobi structure can be obtained out of a Poisson structure (called ``Poissonization'') by the above procedure.
Now, it is well-known that the cotangent $T^* X $ of a Poisson manifold $(X,{\mathcal P})$ is endowed with a Lie natural algebroid bracket with anchor ${\mathcal P}^\# :T^*X \to TX $, whose restriction to exact one-forms is given by $ \lb{\diff f, \diff g}_{\mathcal P}= \diff \, {\mathcal P}[f,g] $. Also, the pair $(TX,T^*X)$ is what is called a Lie bialgebroid over $X$ (one can consult, for instance, \cite{MK} for a detailed introduction to these matters). A natural question, first addressed in \cite{Ker}, is to figure out whether a similar construction can be done for Jacobi structures on $M$. The answer goes as follows: the restriction of $T^* X_{|_M} \to M$ to $M$ can again be endowed with a Lie algebroid structure (which is not a Lie subalgebroid of $T^* X \to X$), that
 forms a Jacobi bialgebroid structure, when paired with some natural algebroid structure that appears on the restriction to $M$ of $TX$. All these constructions can be extended from Poisson structures on Poisson manifolds to Poisson structures on Lie algebroids without any difficulty.

In the continuation of \cite{Kosmann90,Roytenberg}, and in the same spirit of unification and simplification, the purpose of the present short note is to redo, with the help of the big bracket, the theory of Jacobi manifolds and Jacobi (bi)-algebroids, exposed mainly in two articles, \cite{IM} by Iglesias and Marrero, \cite{GM} by Grabowski and Marmo, and continued in \cite{AJR} by Caseiro, de Nicola and Nunes da Costa. The big bracket is the canonical Poisson structure on the graded cotangent of a given vector bundle (considered itself as a graded manifold, see below). It is an unifying tool of a remarkable efficiency when it comes to Lie algebroids or one of the many closely related objects (i.e., for instance, bialgebroids, quasi-Lie bialgebroids, triangular Lie bialgebroids, Courant algebroids, to mention a few), yielding both a nice interpretation of their definitions and mechanical proofs of the theorems that they satisfy. More precisely, our claim is that the big bracket reduces to mechanical (but somewhat cumbersome) computations most known results about Jacobi manifolds, once the objects that these articles deal with have been translated in terms of supermanifolds, a translation that we present herein.
  By mechanical, we mean that it reduces proofs to a play of substitutions involving the Jacobi and Leibnitz identity of the big bracket: we do not claim, however, that the big bracket was the adequate context to guess them.

\bigskip
\noindent
{\bf Acknowledgment.} This work was partially supported by grant no. PTDC/MAT/ 099880/2008 of the Funda\c{c}\~{a}o para a Ci\^encia e a Tecnologia FCT. We would like to thank Joana Nunes da Costa and David Iglesias-Ponte for useful discussions and comments.
\section{The big bracket}

\subsection{Definitions}
\label{sec:defi}

Let $A \to M$ be a vector bundle. There is a natural structure of graded $2$-manifold on $T^* \Pi A$
(i.e. the cotangent space of the supermanifold $\Pi A$) \cite{Roytenberg,Tulczyjew}. Moreover, the (sheaf of) algebra of functions $\superF{A}:=\FF(T^* \Pi A)$  is equipped with a graded Poisson structure which is
called the big bracket and denoted by $\bb{\hbox{$\cdot$,$\cdot$}}$.

We do not intend to give the complete construction of this Poisson algebra, and refer to \cite{Roytenberg} for a more involved introduction, but we recall a few facts on these structures. First, the (sheaf of) algebra $\superF{A}$ admits a $(\N \times \N)$-valued bi-degree:
$$   \superF{A}=\oplus_{k,l \in \N \times \N} \superF{A}^{k,l} $$
 and is graded commutative w.r.t. the total degree, i.e. $ F_1 F_2 = (-1)^{ (k_1+l_1)(k_2+l_2)}F_2 F_1$
 for every $F_1 \in \superF{A}^{k_1,l_1} $ and  $F_2 \in \superF{A}^{k_2,l_2}$.
 The big bracket is a (local) bilinear map $  \superF{A} \times  \superF{A} \to  \superF{A} $,
 denoted $(F_1,F_2) \mapsto \bb{F_1,F_2}$, mapping
 $\superF{A}^{k_1,l_1} \times \superF{A}^{k_2,l_2}  $ to $ \superF{A}^{k_1+k_2-1,l_1+l_2-1}$
 for all $k_1,k_2,l_1,l_2 \in \N$, (with the understanding that $\superF{A}^{k,l}=0$ if $k<0$ or $l<0$)
and which satisfies the following relations
for every $F_i \in \superF{A}^{k_i,l_i}$, $i=1,2,3$,
 \begin{eqnarray}
  \bb{F_1,F_2}&=& -  (-1)^{(k_1+l_1)(k_2+ l_2) } \bb{F_2,F_1} \quad \hbox{(graded skew-symmetry)} \nonumber \\
\bb{F_1, F_2F_3}&=& \bb{F_1,F_2}F_3+(-1)^{(k_1+l_1)(k_2+l_2)}F_2\bb{F_1,F_3} \quad \hbox{(Leibnitz rule)}  \label{eq:Leibnitz} \\
\bb{F_1,\bb{F_2,F_3}}&=& \bb{\bb{F_1,F_2},F_3}+(-1)^{(k_1+l_1)(k_2+l_2)}\bb{F_2,\bb{F_1,F_3}}
 \quad \hbox{(Jacobi)} \label{eq:Jacobi}
 \end{eqnarray}
Explicitly, upon fixing local coordinates $x_i,p^i,\xi_j,\theta^j$ (with $i=1,\dots,n $ and
$j=1,\dots,d $), of respective bi-degrees $(0,0),(1,1),(0,1)$ and $(1,0)$, the
algebra of functions $\superF{A}$ is the graded commutative algebra in those variables, admitting a polynomial dependence in the even variable $p_1,\dots,p_n$ (by skew-symmetry, the dependence in the variables $\xi_1,\theta^1, \dots, \xi_d,\theta^d$ is also polynomial). The big bracket is given in coordinates by:
 $$ \bb{p^i,x_i}=\bb{\theta^j,\xi_j}=1  \quad  i =1, \dots, n, \, \, j=1, \dots , d $$
while all the remaining brackets of coordinate functions vanish.

In this article, we shall use mainly the following points:
\begin{enumerate}
\item There is a natural identification between the algebras $\superF{A}^{0,0}$ and $\FF(M)$.
\item There is a natural identification of graded algebra between $\sum_{k \in \N} \superF{A}^{k,0}$ and $\sum_{k \in \N} \Gamma(\wedge^k A)$. The restriction of $\bb{.,.}$ to this subalgebra is trivial.
\item There is a natural identification of graded algebra between $\sum_{k \in \N} \superF{A}^{0,k}$ and $\sum_{k \in \N} \Gamma(\wedge^k A^*)$. The restriction of $\bb{.,.}$ to this subalgebra is trivial.
\item There is therefore a natural inclusion of $\Gamma(\bigwedge (A \oplus A^*)) \simeq \Gamma(\bigwedge A \otimes \bigwedge A^*)  $
in $\superF{A} $. From now, this inclusion shall be implicitly done, and no notational distinctions shall be made between an element
in  $\Gamma(\bigwedge (A \oplus A^*))$ and its image in $\superF{A}$.
 \item The big bracket between a section of $A$ and a section of $A^*$ is given by the natural pairing;  in equation  $\bb{\xi,\theta}= \bb{\theta,\xi} = \xi (\theta)$, for every $\theta \in \Gamma (A), \xi \in \Gamma(A^*)$.
\item There is a canonical isomorphism of graded Poisson algebras $ \Phi: \superF{A}  \simeq \superF{A^*} $
 intertwining $\superF{A}^{k,l} $ and  $\superF{A^*}^{l,k} $ for all $k,l \in \N$,
\item Let $ F \in \superF{A}^{k,l}$, with $k \geq 1$. If $\bb{\bb{\bb{F,a_1}, \dots }, a_k }=0 $ for all $a_1, \dots, a_k \in \Gamma(A)$, then $F=0$.
\end{enumerate}

There exists an unique function in $\superF{A}$, that we shall denote $\id_A $, which
corresponds to the identity map of the vector bundle $A$, seen as an element of $\Gamma(A^* \otimes A) \subset \superF{A}$. It is explicitly defined by the global section of $\Gamma(A^* \otimes A)  $ given by
\begin{equation}\label{eq:identities}  \id_A := \sum_{j=1}^d \, \xi_j \, \theta^j ,\end{equation}
where $ \xi_1, \dots,\xi_d$ and $\theta^1, \dots , \theta^d  $ are local basis of $A^*$ and $A$
dual to each other.
Taking the bracket with $\id_A$ is a manner to count the bi-degree, more precisely,
for every $P \in \Gamma(\wedge^k A)$, $\Psi \in \Gamma (\wedge^l A^*)$,
\begin{equation}\label{eq:ligacao com degree} \bb{P \Psi, \id_A} :=  (k-l) P \Psi .\end{equation}
Of course, one can also consider the function $\id_{A^*} \in \superF{A^*}$.
Under the canonical isomorphism $\Phi:\superF{A}\simeq\superF{A^*} $, both functions are related by:
\begin{equation}\label{eq:ligacao com dualidade}  \Phi (\id_A) = - \id_{A^*}.\end{equation}


\subsection{Preliminary results}
\label{sec:prelim}

For a future purpose, we wish to establish several facts about the behaviour of the big bracket when one adds a copy of $\R$ either to the base or to the fibers.

\smallskip
\noindent
{\bf A. Enlarging the base.} Let $A \to M$ be a vector bundle. Denote by $p^* A \to (M \times \R)$ the pull-back of $A \to M$ through the projection onto the first component $p: M \times \R \to M$. We denote by $t \in \FF (M \times \R) \subset \superF{p^* A}$ the projection onto the second component. There is a canonical inclusion $ \mathfrak{i} : \superF{A} \hookrightarrow \superF{ p^* A} $.
  There exists an unique function $\partialt \in \superF{p^* A}$ such that
\begin{equation}\label{eq:def_partialt}
\bb{\partialt, \mathfrak{i} (F) }=0  \hbox{ for all $F \in \superF{A}$ and }  \bb{\partialt,t}=1 .
\end{equation}
The bi-degree of this function is $(1,1)$. With a slight abuse of notations, we shall consider $\superF{A}$ as a subset of $\superF{p^* A}$, erasing therefore the canonical inclusion $ \mathfrak{i} $.

Notice that $\partialt$ is an even function, and that, for every $F_0,\dots,F_k\in \superF{A}$, the relation $\sum_{i=0}^k F_i ({\partialt})^i =0$ holds if and only if $F_0=\dots=F_k=0$.

\smallskip
\noindent
{\bf B. Enlarging the fibers.} Let $B \to M $ be a vector bundle. We call $A=B \oplus \R \to M $ the direct sum of $B$  with the trivial bundle $\R \times M \to M$. There is a natural inclusion $ \mathfrak{j}:\superF{B} \subset \superF{A}$, which preserves the big bracket.
 We define $\phi \in  \Gamma(A^*) \subset \superF{A}$ to be the projection onto the second component, i.e the section of $A^*= (B \oplus \R)^*$ defined by $\phi ( b,f)=f $ for all $b \in \Gamma(B), f \in \FF (M)$. The following is immediate:
 \begin{equation}\label{eq:prop_phi}
\bb{\phi, \mathfrak{j} (F) }=0  \hbox{ for all $F \in \superF{B}$ and }  \bb{\phi , \epsilon}=1 ,
\end{equation}
 where $\epsilon$ is the section of $A:= B \oplus \R \to M$ given by $m \to (0_m,1) \in A_m = B_m \oplus \R$ (in the previous $0_m$ is the zero element in the vector space $B_m$). With a slight abuse of notations, we shall consider $\superF{B}$ as a subset of $\superF{A}$, erasing therefore the canonical inclusion $ \mathfrak{j} $.


\section{Lie and Jacobi algebroids, Lie bialgebroids and Jacobi bialgebroids.}
\label{sec:Lie_e_Jac}

\subsection{ (Pre-)Jacobi algebroids.}

\label{sec:pre-jac}
Lie algebroids and pre-Lie algebroids are in general introduced through brackets and anchors, however,
it is well-known that the supergeometric point of view is strictly equivalent to those more classical ones,
see \cite{Roytenberg}.

\begin{defi}\label{def:algebroids}
Let $A \to M$ be a vector bundle. A \emph{pre-Lie algebroid (on $A$)} is a function in $ \superF{A}^{1,2}$
(i.e. a function of bi-degree $(1,2)$).

The \emph{bracket of a pre-Lie algebroid $\mu$} is the bilinear endomorphism of $ \Gamma(\bigwedge A):=\oplus_{k \in \N} \Gamma(\wedge^k A)$ defined, for all $P,Q \in  \Gamma(\bigwedge A)$ by
\begin{equation}\label{pseudoGerstenhaber} \lb{P,Q}_\mu= \bb{\bb{P,\mu},Q} .\end{equation}
The \emph{differential} of a pre-Lie algebroid $\mu$ is the linear endomorphism of $\Gamma(\bigwedge A^*):=\oplus_{k \in \N} \Gamma(\wedge^k A^*) $, of degree $+1$, defined, for all $\Psi \in \Gamma(\bigwedge A^*)$, by $\diff_\mu (\Psi):= \bb{ \mu, \Psi } $.
\end{defi}

Recall that the restriction to $\Gamma(A) \times \FF(M) \to \FF(M)$ of the bracket $\LB_\mu$ is of the form
 $$ (a,f) \mapsto \rho_\mu (a) \, [f]  \quad \forall a \in \Gamma (A), f \in \FF (M),$$
 for some  vector bundle morphism $ \rho_\mu : A \to TM $ (over the identity of $M$)
 called the \emph{anchor map}.
 In the previous, the notation $f \mapsto X[f]$ is used to denote the derivation of $\FF(M)$ associated to a vector field $X$ on $M$ (a notation that we shall also use for multi-vector fields).

 A pre-Lie algebroid $\mu$ (on $A$) is said to be a \emph{Lie algebroid (on $A$)} if $\bb{\mu,\mu}=0$.
   It is a classical result (see, for instance, \cite{Kosmann90,Roytenberg}) that a pre-Lie Algebroid is a Lie algebroid if and only if one of the following equivalent conditions is satisfied: \emph{(i)}  (\ref{pseudoGerstenhaber}) satisfies the Jacobi identity (hence is a Gerstenhaber bracket) \emph{(ii)} $\diff_\mu^2=0$ or,  \emph{(iii)}, the restriction of $\LB_\mu$ to $\Gamma(A) \times \Gamma(A) \to \Gamma (A)$ is a Lie algebra bracket.

   We now define pre-Jacobi algebroids.

\begin{defi}\label{def:jacobi_algebroids}
Let $A \to M$ be a vector bundle. A \emph{pre-Jacobi algebroid (on $A$)} is a function of $ \superF{A}$ of the form $ \mu+\phi$, with $\mu \in \superF{A}^{1,2}$  and $\phi \in \superF{A}^{0,1} $ (said differently $\mu $ is a pre-Lie algebroid while $\phi$ is a section of $A^*$).

The \emph{bracket of a pre-Jacobi algebroid $\mu$} is the bilinear endomorphism of $ \Gamma(\bigwedge A):=\oplus_{k \in \N} \Gamma(\wedge^k A)$ defined, for all $P,Q \in  \Gamma(\bigwedge A)$ by
 $$ \lb{P,Q}_{\mu,\phi} := \bb{ \bb{ P, \mu + \id_{A} \phi }, Q} - P \bb{\phi, Q}+ \bb{P,\phi} Q $$
 (see section \ref{sec:defi} for the definition of $\id_A$).
The \emph{differential} of a pre-Jacobi algebroid $\mu+\phi$ is the linear endomorphism $\diff_{\mu,\phi} $ of $\Gamma(\bigwedge A^*):=\oplus_{k \in \N} \Gamma(\wedge^k A^*) $
defined, for all $\Psi \in \Gamma(\bigwedge A^*)$, by $\diff_{\mu,\phi} (\Psi):= \bb{ \mu, \Psi } + \phi \Psi $.
\end{defi}

\begin{lem}\label{lem:explicit}
The bracket of a pre-Jacobi algebroid is given explicitly, by the following formula,
for all $P \in \Gamma(\wedge^k A), Q \in \Gamma(\wedge^l A)$:
 \begin{equation}\label{ex:explicit} [P,Q]_{\mu,\phi} = \lb{P,Q}_\mu  + (k-1) P (\imath_\phi Q ) +(-1)^k (l-1) (\imath_\phi P ) Q \end{equation}
\end{lem}
\begin{pf}
Using the definitions of the pre-Lie brackets $[\hbox{$\cdot$,$\cdot$}]_{\mu,\phi}$ and
$[\hbox{$\cdot$,$\cdot$}]_{\mu}$, we obtain:
\begin{eqnarray*}    [P,Q]_{\mu,\phi} &=& \bb{ \bb{ P, \mu }, Q }+\bb{ \bb{ P, \id_A \phi }, Q }
- P \bb{\phi, Q}+ \bb{P,\phi} Q  \\
 &=&  \bb{ \bb{ P, \mu }, Q }+\bb{ \bb{ P, \id_A \phi }, Q }
- P \imath_\phi  Q- (-1)^k(\imath_\phi P) Q\\
 &=&  \lb{P,Q }_\mu+\bb{ \bb{ P, \id_A \phi }, Q }
- P \imath_\phi  Q- (-1)^k(\imath_\phi P) Q.
\end{eqnarray*}
By making several use of the Leibnitz identity (\ref{eq:Leibnitz}), one computes
  \begin{eqnarray*}\bb{ \bb{ P, \id_A \phi }, Q }&=& \bb{\bb{P,\id_A }\phi, Q} +  \bb{ \bb{P,  \phi }\id_A ,Q} \\
   & = &      k \bb{P \phi , Q} +\bb{ \bb{P,  \phi } \id_A ,Q} \quad \hbox{  by (\ref{eq:ligacao com degree})} \\
   & = &    k P\bb{\phi , Q} +\bb{ \bb{P,  \phi } \id_A ,Q}\\
   &=&   k P\bb{\phi , Q} +\bb{P,  \phi } \bb{\id_A ,Q} \\
   & = &  k P \bb{\phi , Q} - l \bb{P,  \phi } Q \quad \hbox{  by (\ref{eq:ligacao com degree})} \\
    & = &  k P \imath_\phi Q + l (-1)^k(\imath_\phi P) Q    \end{eqnarray*}
This completes the computation.
\end{pf}

\begin{rmk}In Theorem 3.5 \cite{IM} or Equation (25) \cite{GM}, a bracket is constructed out of the datae defining a Jacobi algebroid. Comparing their quite explicit formulas with formula (\ref{ex:explicit}) above proves the coincidence of our bracket with theirs (more precisely, the match is exact with \cite{GM}, but is only up to signs with \cite{IM} where the Gerstenhaber bracket of a Lie algebroid is given by $ P,Q \mapsto \bb{\bb{\mu,P},Q}$ and not by (\ref{pseudoGerstenhaber})).
\end{rmk}

 We also easily deduce  the following properties from those of the big bracket, for all
 $a \in \Gamma(\wedge^k A), b \in \Gamma(\wedge^l A),c \in \Gamma(\wedge^m A)$:
   $$ \begin{array}{rcl} \lb{b,a}_{\mu,\phi} &=& -(-1)^{kl} \lb{a,b}_{\mu,\phi}\\
 \lb{a,b c}_{\mu,\phi}  &=&  \lb{a,b}_{\mu,\phi}  c+(-1)^{lm}  \lb{a, c}_{\mu,\phi} b - \bb{\phi,a} bc \end{array}.$$

\begin{defi}
 A pre-Jacobi algebroid $\mu+\phi$ (on $A$) is said to be a \emph{Jacobi algebroid (on $A$)} when $\bb{\mu+\phi,\mu+\phi}=0$.
\end{defi}

Spelled out, the condition $\bb{\mu+\phi, \mu+\phi}=0 $ yields the two conditions $ \bb{\mu,\mu}=0$ and $\bb{\mu,\phi}=\diff_\mu \phi=0$. The first of these conditions means that $\mu$ is a Lie algebroid, and the second one means that $\phi$ is a cocycle of this Lie algebroid, i.e. that $\phi(\lb{a,b}_\mu) = \rho_\mu (a) \, [ \phi(b)]- \rho_\mu (b) \, [\phi(a)]$ for all $a,b \in \Gamma(A)$ ($\rho_\mu: A \to TM$ being the anchor map defined above). In conclusion, a Jacobi algebroid is a Lie algebroid endowed with an algebroid $1$-cocycle, which is the usual definition (compare with \cite{IM}, where such an object is called
``Lie algebroid in the presence of a $1$-cocycle'').

 The next result appeared in both \cite{IM} and \cite{GM}:
we prove it here with the help of super-geometric formalism.

\begin{prop}
\label{prop:Gra}
Let $A\to M$ be a vector bundle. For every pre-Jacobi algebroid $\mu+\phi$ on $A$, the following are equivalent:
\begin{enumerate}
\item[(i)] $\mu+\phi$ is a Jacobi algebroid structure on $A$;
\item[(ii)] the operator $\diff_{\mu,\phi}$ squares to $0$;
\item[(iii)] $\LB_{\mu,\phi}$ satisfies the graded Jacobi identity:
  $$ (-1)^{km}\lb{\lb{a,b}_{\mu,\phi},c}_{\mu,\phi} + (-1)^{lk}\lb{\lb{b,c}_{\mu,\phi},a}_{\mu,\phi} +(-1)^{ml}\lb{\lb{c,a}_{\mu,\phi},b}_{\mu,\phi}=0  ,$$
for all homogeneous $a,b,c \in \Gamma(\bigwedge A)$ of degrees $k,l,m$ respectively.
\end{enumerate}
\end{prop}
\begin{pf}
The equivalence between \emph{(i)} and \emph{(ii)} follows from the relation $  \diff_{\mu,\phi}^2 (\Psi)= \frac{-1}{2}\bb{\bb{\mu+\phi,\mu+\phi}  ,\Phi}$, together with the seventh result listed in section \ref{sec:defi}.
 Let us prove the equivalence between \emph{(i)} and \emph{(iii)}. A direct computation with the help of the usual properties of the big bracket gives the following expression for the Jacobiator the bracket $\LB_{\mu,\phi}$:
\begin{eqnarray*}(-1)^{lm} \lb{\lb{a,b}_{\mu,\phi},c}_{\mu,\phi} + \hbox{ c.p. } &=& (-1)^{lm}\lb{\lb{a ,b}_\mu ,c}_\mu +\hbox{ c.p. }\\ & & \quad +
(-1)^{lm}(m-1)( \imath_\phi \lb{ a, b}_\mu -  \lb{\imath_\phi a, b}_\mu- \lb{ a,\imath_\phi b}_y  ) c +\hbox{ c.p. } \end{eqnarray*}
(where $+\hbox{ c.p. }$ indicates that we add all the terms obtained by permuting the variables $a,b,c$).
Now, the Jacobi identity of the big bracket implies that:
\begin{eqnarray*} \imath_\phi \lb{ a, b}_\mu -  \lb{\imath_\phi a, b}_\mu- \lb{ a,\imath_\phi b}_\mu
& = &  \imath_\phi \bb{ \bb{a,\mu}, b} -  \bb{\bb{ \imath_\phi a,\mu} b}- \bb{ \bb{ a,\mu} \imath_\phi b}  \\
& = &  \bb{\phi, \bb{ \bb{a,\mu}, b}} -  \bb{\bb{ \bb{\phi ,a},\mu}, b}_\mu- \bb{ \bb{ a,\mu} ,\bb{\phi, b}}  \\
& = &  \bb{ \bb{a,\bb{\phi,\mu}}, b}
\end{eqnarray*}
As a consequence:
$$(-1)^{lm} \lb{\lb{a,b}_{\mu,\phi},c}_{\mu,\phi} + c.p. = (-1)^{lm}\lb{\lb{a ,b}_\mu ,c}_\mu +
(-1)^{lm}(m-1)  \bb{ \bb{a,\bb{\phi,\mu}}, b}  c +\hbox{ c.p. }$$
If $\mu$ is a Lie algebroid, the first term on the right hand side vanishes, and if $\bb{\phi,\mu}=0$, the second term vanishes as well, so that the graded Jacobi identity is satisfied. Conversely, if the graded Jacobi identity is satisfied
for all homogeneous $a,b,c$, then, choosing $c=1$, it follows from the previous identity that
 $$ \bb{ \bb{a,\bb{\phi,\mu}}, b} =0 $$
 for all homogeneous $a,b$. Since $\bb{\phi,\mu} \in \superF{A}^{0,2}$, this implies  $\bb{\phi,\mu}=0$
 (see the seventh result listed in section \ref{sec:defi}). In turn, this implies that
 the bracket $\LB_\mu$ satisfies the graded Jacobi identity, a property that holds true if and only if $\mu$ is a Lie algebroid.
\end{pf}

Given a pre-Jacobi algebroid  $\mu+\phi $ on $A$, we define a family indexed by a parameter $c \in \R$
of pre-Lie algebroid structures, denoted $\mu_\phi^{c}  $, by
$$ \mu_\phi^c := e^{-ct}(\mu+ \phi (\partialt + c\id_A  )) $$
where we use the notations of section \ref{sec:prelim}-A, i.e.:
\begin{enumerate}
\item $p^* A \to (M \times \R)$ is the pull-back of $A \to M$ through the projection onto the first component $( M \times \R ) \to M$,
\item $t \in \FF (M \times \R) \subset \superF{p^* A} $ is the parameter on $\R$ (i.e. the projection onto the second component),
\item $ \superF{A} $ is considered as a subalgebra of $ \superF{p^* A}$,
\item $\partialt \in \superF{A}^{1,1}$ is defined as in section \ref{sec:prelim}-A.
\end{enumerate}
 The anchors and brackets of these structures correspond, for $c=0,1$, to (4.16-4.19) in \cite{IM}. The pre-Lie algebroid $ \mu_\phi^{1} $  (on $p^* A$) is called the \emph{Poissonization} of the pre-Jacobi algebroid  $\mu+\phi$ (on $A$). Notice also that $\mu_\phi^0= \mu +\phi \partialt $.

\begin{lem}\label{lem:exp}
Define $\ad_{t \,\id_A}$ to be the linear endomorphism of $\superF{A} $ defined by $F \mapsto \bb{ t \, \id_A , F} $. For all $c,x \in \R$, the following relation holds:
   $$ {\rm  exp} \, (\ad_{-xt \, \id_A}) \mu^c_\phi= \mu^{c+x}_\phi .   $$
\end{lem}
\begin{pf}
A direct computation gives
 $$ \ad_{t \,\id_A} \mu^{0}_\phi = t \mu_\phi^0 - \phi \id_A  \mbox{ and }  \ad_{t \, \id_A} \phi \id_A = t \phi  \id_A.$$
Now, it is a general fact that for every vector space $E$, every $L \in {\rm End} (E,E)$ and every
${\mathfrak a}, {\mathfrak b}\in E$, if $L({\mathfrak a})= t{\mathfrak a} - {\mathfrak b}$ and $L({\mathfrak b})=t{\mathfrak b} $, we have ${\rm exp} \, (-cL) \, {\mathfrak a} = e^{-tc} {\mathfrak a} + c e^{-tc} {\mathfrak b} $.
Applied to $E:= \superF{A} $, $L:= \ad_{t \, \id_A} $, ${\mathfrak a}:= \mu^0_\phi$, ${\mathfrak b}:=\phi \id_A $, this gives
 $$  {\rm  exp} \, (\ad_{-ct \,\id_A}) \mu^0_\phi= e^{-ct}\mu_\phi^0 + c e^{-ct} \id_A \phi = \mu_\phi^c . $$
Applying $ {\rm  exp} \, (\ad_{-xt \,\id_A})$ to both sides of this relation gives:
  $$ {\rm  exp} \, (\ad_{-xt \,\id_A}) \, {\rm  exp} \, (\ad_{-ct \,\id_A}) \mu^0_\phi=  {\rm  exp} \, (\ad_{-xt \,\id_A}) \, \mu_\phi^c ,$$
  The relation
  $${\rm  exp} \, (\ad_{-xt \,\id_A}) \, {\rm  exp} \, (\ad_{-ct \,\id_A}) \,  \mu^0_\phi= {\rm  exp} \, (\ad_{-(x+c)t \,\id_A}) \, \mu^0_\phi = \mu^{c+x}_\phi,$$
  now gives the required result.
\end{pf}

The next lemma is dealt with in \cite{IM} for the cases $c=0,1$.

\begin{lem}\label{lem:alg_associadas_Jaco_alg}
Let $A \to M$ be a vector bundle,  $\mu+\phi$ a pre-Jacobi algebroid on $A\to M$. The following are equivalent ($p^* A, \mu_{ \phi}^{c}$ being as defined above):
\begin{enumerate}
\item[(i)] $\mu+\phi$ is a Jacobi algebroid on $A $;
\item[(ii)] There exists $c \in \R$ such that the pre-Lie algebroid $ \mu_\phi^{c}$ is a Lie algebroid on $p^* A $;
\item[(iii)] For all $c \in \R$, the pre-Lie algebroid $\mu_\phi^{c}$ is a Lie algebroid on $ p^* A$.
\end{enumerate}
\end{lem}
\begin{pf}
 Using the usual properties of the big bracket, one computes
$$ \begin{array}{ll}\bb{ \mu_{ \phi}^{0}, \mu_{ \phi}^{0}} &=  \bb{\mu + \phi \partialt, \mu + \phi \partialt}
   \\  &= \bb{\mu,\mu}+ 2 \bb{\mu,\phi}\partialt
     \end{array}
     $$
The vanishing of $\bb{ \mu_{ \phi}^{0}, \mu_{ \phi}^{0}} $ is therefore
(in view of section \ref{sec:prelim}-A) tantamount to the vanishing of both $\bb{\mu,\mu}$ and $\bb{\mu,\phi}$. Hence $\mu^0_\phi  $ is a Lie algebroid if and only if $\mu +\phi$ is a
 Jacobi algebroid. Now, in view of lemma \ref{lem:exp}, we  have for all $c \in \R$
 $$\bb{\mu_\phi^c,\mu_\phi^c } =  \bb{{\rm  exp} \,(\ad_{-ct\, \id_A}) \mu_\phi^0 ,{\rm  exp}\, (\ad_{-ct\, \id_A}) \, \mu_\phi^0 }  =  {\rm  exp} \,(\ad_{-ct \,\id_A}) \, \bb{ \mu_\phi^0, \mu_\phi^0}  $$
In particular, $ \mu_\phi^c $ is a Lie algebroid if and only if $\mu_\phi^0$
is a Lie algebroid. This completes the proof.
\end{pf}

\subsection{Jacobi bialgebroids.}

Upon identifying $\superF{A}$ and $\superF{A^*} $ with the help of $\Phi$, it is possible (since $\Phi$ intertwines  $\superF{A}^{k,l} $ and $ \superF{A^*}^{l,k} $ for all $k,l \in \N$, see  section \ref{sec:defi}) to consider pre-Lie algebroids structures on $A^*$ as elements of  $\superF{A}^{2,1} $ and pre-Jacobi algebroids structure on $A^*$ as functions in $\superF{A}$ of the form $\nu+X$, with $\nu \in \superF{A}^{2,1} $ and $X \in  \superF{A}^{1,0}$.

Now, recall \cite{Kosmann90} that a Lie bialgebroid is a pair of pre-Lie algebroid structures $\mu,\nu$ on $A\to M$ and on $A^* \to M$ respectively (i.e. functions $\mu \in \superF{A}^{2,1}$ and $\nu \in \superF{A}^{1,2}$ respectively), such that $\bb{\nu+\mu,\nu+\mu}=0 $. Requiring this condition is equivalent to require that both pre-Lie algebroid structures are indeed Lie algebroids, and that the following compatibility condition is satisfied:
 $$ \diff_{\nu} \lb{ a,b }_\mu =  \lb{  \diff_{\nu}a,b }_\mu +(-1)^{k-1}\lb{ a,\diff_{\nu}  b }_\mu \quad \forall a \in \Gamma(\wedge^k A), b \in \Gamma(\wedge A)  $$
where $\diff_{\nu}$, $\LB_\mu $ are as in definition \ref{def:algebroids}.
Following the same idea, Jacobi bialgebroids are introduced in \cite{GM} as follows.

\begin{defi}
Let $A \to M$ be a vector bundle. A \emph{pre-Jacobi bialgebroid} is a pair $(\mu+\phi,\nu+X) $, where
\begin{enumerate}
\item $\mu+\phi$ is a pre-Jacobi algebroid structure on $A$,
\item $\nu+X$ is a pre-Jacobi algebroid structure on $A^*$,
\end{enumerate}
such that the following compatibility condition is satisfied
    \begin{equation}\label{cond_Gra}
    \diff_{\nu,X} \lb{ a,b }_{\mu,\phi} =  \lb{  \diff_{\nu,X}a,b }_{\mu,\phi} +(-1)^{k-1}\lb{ a,\diff_{\nu,X}  b }_{\mu,\phi}  \quad \forall a \in \Gamma(\wedge^k A), b \in \Gamma(\wedge A) ,\end{equation}
where $\LB_{\mu,\phi},\diff_{\nu,X}$ are as in definition \ref{def:jacobi_algebroids}.
A \emph{Jacobi bialgebroid} is a pre-Jacobi bialgebroid $(\mu+\phi,\nu+X)$ such that
$\mu+\phi$ is a Jacobi algebroid and $ \diff_{\nu,X}$ squares to $0$.
\end{defi}

We can now give the main result of this section, which is a characterization, in terms of the big bracket,
of Jacobi bialgebroid structures.

\begin{theo}\label{theo:main}
Let $A \to M $ be a vector bundle, $\mu+\phi$ a pre-Jacobi algebroid structure on $A$ and $\nu+X$ a pre-Jacobi algebroid structure on $A^*$. Set
  \begin{equation}\label{eq:horrible}\left\{\begin{array}{rcl} \diamondsuit &:= &\bb{\phi,X}= \phi(X) \\ \clubsuit &:= & \bb{\mu,X} + \bb{\nu,\phi} \\ \spadesuit & := & \bb{\mu,\nu}+ \mu X+\nu\phi - \id_{A}\left(  -\phi X + \frac{\bb{ \mu,X }- \bb{ \nu,\phi }}{2}\right). \\
 \end{array}\right.\end{equation}
 \begin{enumerate}
 \item
The pair $(\mu+\phi,\nu+X) $ is a pre-Jacobi bialgebroid if and only if
$$ \diamondsuit=\clubsuit=\spadesuit=0 .$$
\item The pair $(\mu+\phi,\nu+X) $ is a Jacobi bialgebroid if and only if
$$ \diamondsuit=\clubsuit=\spadesuit= \bb{\mu+\phi,\mu+\phi} =  \bb{\nu+X,\nu+X}=0 .$$
\end{enumerate}
\end{theo}
\begin{pf}
We prove the first item. The idea is to express, with the help of the big bracket, the quantity
$$ H_{\lfloor .,.\rfloor, D } (a,b):= D \lfloor a,b \rfloor -  \lfloor Da,b \rfloor -(-1)^{k-1} \lfloor a, D b \rfloor , \quad  a \in \Gamma(\wedge^k A), b \in \Gamma( \bigwedge A) $$
where $D$ stands either for $\ad_\nu:= \bigwedge A \mapsto \bigwedge A$ or for the left-multiplication by $X$, i.e. $m_X(b):= Xb$ , and $ \lfloor a,b \rfloor $ stands either for the derived bracket
$ \lb{a,b}_{\mu'}:=\bb{\bb{a, \mu'},b } $, with $\mu':= \mu+\id_A \phi $, or for the assignment $ (( a,b))_\phi := \bb{a,\phi}b -a\bb{\phi,b}$.
In fact, by definition, $ \diff_{\nu,X}:=\ad_\nu + m_X $, while $ \LB_{\mu,\phi}  := \LB_{\mu} +(( .,.))_\phi $, hence
   \begin{equation}\label{eq:desab1}
   \begin{array}{ll} &\diff_{A^*,\phi} \lb{ a,b }_{A,\phi} -  \lb{  \diff_{A^*,\phi}a,b }_{A,\phi}  -(-1)^{k-1}\lb{ a,\diff_{A^*,\phi}  b }_{A,\phi} \\=& \big(  H_{\LB_{\mu'},\ad_\nu} +H_{\LB_{\mu'}, m_X } +H_{ (( .,.))_\phi , \ad_\nu } +H_{ (( .,.))_\phi ,m_X  }\big) (a,b)
   \end{array}\end{equation}
A cumbersome but straightforward computation, involving only the Leibnitz and Jacobi identity of the big bracket
(and valid for arbitrary $\mu' \in \superF{A}^{1,2},\nu \in  \superF{A}^{2,1}, X  \in  \superF{A}^{1,0},\phi \in  \superF{A}^{0,1} $), yields to the following results:
$$
\left\{ \begin{array}{rcl}
  H_{\LB_{\mu'},\ad_\nu} (a,b)        & = & {\mathcal D}(a,(-1)^k\bb{\nu,\mu'},b) \\
  H_{\LB_{\mu'}, m_X } (a,b)          & = & {\mathcal D}(a,(-1)^k  \mu' X,b) +  {\mathcal E}(a,(-1)^k \bb{X,\mu'},b)   \\
  H_{ ((.,.))_\phi , \ad_\nu } (a,b)  & = & {\mathcal E}(a,(-1)^k\bb{\nu,\phi},b) \\
  H_{ (( .,.))_\phi ,m_X  } (a,b)     & = & {\mathcal F}(a,2(-1)^{k+1} \bb{\phi,X} ,b) + {\mathcal E}(a,-(-1)^k X\phi,b ) \\
\end{array} \right.
$$
where, for all $a,b \in \bigwedge A$, for all function $F \in \superF{A}$ of bi-degree $(1,2)$,
$$ {\mathcal D}(a,F,b):= \bb{\bb{a, F},b } \,\, , \,\, {\mathcal E}(a,F,b):= \bb{a, F}b -  a \bb{F,b}  \,\, , \,\,{\mathcal F}(a,F,b):= F ab .$$
Introducing these results in (\ref{eq:desab1}) yields
  \begin{equation}\label{eq:etape2}\begin{array}{rcl}  \diff_{A^*,\phi} \lb{ a,b }_\phi -  \lb{  \diff_{A^*,\phi}a,b }_\phi -(-1)^{k+1}\lb{ a,\diff_{A^*,\phi}  b }_\phi &=& (-1)^k {\mathcal D}(a, \hat{\spadesuit} ,b) \\
& &+ (-1)^k  { \mathcal E}(a, \hat{\clubsuit} ,b)\\
& &+ (-1)^k {\mathcal F}(a, \hat{\diamondsuit} ,b)\\
\end{array}  \end{equation}
where
\begin{equation}\label{eq:etape3}
\left\{ \begin{array}{rclcl}
\hat{\spadesuit}    & := &  \bb{\nu,\mu'}+ \mu' X & = &   \spadesuit  - \frac{1}{2}\id_A \clubsuit   \\
\hat{\clubsuit}     & := &  \bb{X,\mu'} + \{\nu,\phi\} - X \phi &=&   \bb{\mu,X} + \{\nu,\phi\}+ \id_A \bb{X,\phi}= \clubsuit+\id_A \diamondsuit\\
\hat{\diamondsuit}  & := & -2 \bb{\phi,X} & = & -2 \diamondsuit.\\
\end{array} \right.
\end{equation}
From these expressions, it follows that, if $ \diamondsuit=\clubsuit=\spadesuit=0$,
then $\hat{\spadesuit}= \hat{\clubsuit}= \hat{\diamondsuit}=0 $, hence
(\ref{cond_Gra}) holds true, and the structure is a pre-Jacobi bialgebroid structure.

Let us prove the converse. First, notice that ${\mathcal D}(a, \hat{\spadesuit} ,b)=0 $ if $a=1$ or $b=1$,
and ${\mathcal E}(a, \hat{\spadesuit} ,b)=0$ if $ a=b=1$.
Hence, if condition (\ref{cond_Gra}) is satisfied, one sees by plugging $a=b=1$ in (\ref{eq:etape2}) that $  \hat{\diamondsuit}=0$.
Then plugging $a=1$ and letting $b$ be an arbitrary section of $A$, we conclude that
$\bb{\hat{\clubsuit},b} =0$ for all $b \in \Gamma (\bigwedge A)$. Since $\hat{\clubsuit}$ is of bidegree $(1,1)$, according to the seventh item listed in section \ref{sec:defi}, we have $\hat{\clubsuit}=0$. For similar reasons, we deduce from $\bb{\bb{a,\hat{\spadesuit}},b} =0$, valid for $a,b$  arbitrary sections of $A$, that $\hat{\spadesuit}=0$.
In view of (\ref{eq:etape3}), we have $\diamondsuit=\clubsuit=\spadesuit=0 $ which completes the proof.

The second item is an immediate consequence of the first one and of proposition \ref{prop:Gra}
\end{pf}

In the equations (\ref{eq:horrible}), the roles of $A$ and $A^*$ are symmetric, i.e. applying the canonical isomorphism $\Phi: \superF{A} \simeq \superF{A^*}$ (and using (\ref{eq:ligacao com dualidade})), one obtains equations of the same form, which yields the next corollary, which already appears in both \cite{IM} and \cite{GM}.

\begin{coro}\label{cor:simetria}
Let $A \to M $ be a vector bundle, $\mu+\phi$ a pre-Jacobi algebroid structure on $A$, and $\nu+X$ a pre-Jacobi algebroid structure on $A^*$. The pair $(\mu+\phi,\nu+X)$ is a pre-Jacobi bialgebroid structure on $A$ if and only if the pair $(\nu+X,\mu+\phi)$ is a pre-Jacobi bialgebroid structure on $A^*$.
\end{coro}

In lemma \ref{lem:alg_associadas_Jaco_alg} we constructed, out of a Jacobi algebroid structure $\mu+\phi$ on a vector bundle $B \to M$, a family indexed by $c \in \R$ of Lie algebroid structures on the pull-back vector bundle $p^* B \to M \times \R$ ($p: M \times \R \to M$ being the projection on the first component). We denoted these structures by $ \mu_\phi^c$.
 Hence, provided that we give ourself $\mu + \phi$ a pre-Jacobi algebroid structure on a vector bundle $A \to M$ and $ \nu + X $ a pre-Jacobi algebroid structure on the dual bundle $ A^* \to M $, we can construct,
 for every $c,d \in \R$:
 \begin{enumerate}
 \item Lie algebroid structures $\mu_\phi^c$ on $ p^* A \to M \times \R  $,
 \item Lie algebroid structures $\nu_X^d$ on the dual bundle $ p^* A^* \to M \times \R  $.
 \end{enumerate}
Since the vector bundles $p^* A^* \simeq (p^* A)^*$ are dual one to the other,  it is natural to ask whether one can pair these previous structures to form Lie bialgebroids.
The next corollary gives an answer to this question, and generalizes results obtained in \cite{GM,IM}
for $c=0,1$.

\begin{coro}\label{cor:ligacao_com_bialg}
Let $A \to M $ be a vector bundle, $\mu+\phi$ a pre-Jacobi algebroid structure on $A$, and $\nu+X$ a pre-Jacobi algebroid structure on $A^*$. Choose an arbitrary $c \in \R$.
%
\begin{enumerate}
\item
The following are equivalent:
\begin{enumerate}
\item[(i)] the pair $(\mu+\phi,\nu+X  ) $ is a pre-Jacobi bialgebroid on $A$;
\item[(ii)] $\bb{\mu^c_{\phi} , \nu^{1-c}_{X}}=0 $.
\end{enumerate}
\item The following are equivalent:
\begin{enumerate}
\item[(i)] the pair $(\mu+\phi,\nu+X  ) $ is a Jacobi bialgebroid on $A$;
\item[(ii)] the pair $(\mu^c_{\phi} , \nu^{1-c}_{X}) $ is a Lie bialgebroid on $p^* A$.
\end{enumerate}
\end{enumerate}
\end{coro}
\begin{pf}
Item \emph{2)} is a clear consequence of item \emph{1)}, so we only include a proof of the first item.
First:
\begin{eqnarray*}{\rm exp}\,\ad_{-c t\, \id_A }\bb{\mu^c_{\phi} , \nu^{1-c}_{X}} &=&
\bb{{\rm exp}\,\ad_{-c t\, \id_A }\,\mu^c_{\phi} , {\rm{ exp}}\,\ad_{-c t\, \id_A }\,\nu^{1-c}_{X}}\\
&=& \bb{{\rm exp}\,\ad_{-c t\, \id_A }\,\mu^c_{\phi} , {\rm{ exp}}\,\ad_{c t \, \id_{A^*} }\,\nu^{1-c}_{X}}  \\ &=&\bb{\mu_\phi^0,\nu_X^1 } \quad \hbox{ (by lemma \ref{lem:exp})}
\end{eqnarray*}
so that  the vanishing of $ \bb{\mu^c_{\phi} , \nu^{1-c}_{X}}$ is equivalent to the vanishing
of $\bb{\mu^0_{\phi} , \nu^{1}_{X}} =0$.
In view of theorem (\ref{theo:main}), it suffices to prove that the conditions $\clubsuit=\diamondsuit=\spadesuit=0$ are satisfied if and only if
$\bb{\mu^0_{\phi} , \nu^{1}_{X}} =0$. By a direct computation, we obtain:
$$
\begin{array}{rcl}
\bb{\mu_\phi^0,\nu_{X}^1}  &=& \bb{ \mu +  {\partialt} \phi , e^{-t} \big(\nu +( \partialt + \id_{A^*}) X \big) } \\
 & = &  e^{-t}  \bb{\phi,X} \, (\partialt)^2  \\
 & & +  e^{-t}  \partialt \, ( \bb{ \mu,X }+\bb{\phi,\nu}+ \id_{A^*} \bb{\phi,X}  ) \\
 & & +e^{-t}  ( \bb{\mu,\nu}+ \mu X - \phi \nu + \id_{A^*} ( \bb{\mu,X} - \phi X ) )\\
 &= & e^{-t} \big( (\partialt)^2\diamondsuit  + \partialt (\clubsuit +\id_{A^*} \diamondsuit )
 +  \spadesuit - \frac{1}{2}\id_{A^*}\clubsuit \big)
\end{array}
$$
Now, a function $F\in \superF{ p^* A} $ of the form $F= (\partialt)^2 F_1+ \partialt  F_2 + F_3 $, with $F_1,F_2,F_3 \in \superF{A}$ is zero if an only if $F_1=F_2=F_3=0$ (as stated in section \ref{sec:prelim}-A). We can therefore conclude that $\bb{ \mu_\phi^0 , \nu_X^1 } =0$ if and only if $\diamondsuit =\clubsuit +\id_{A^*} \diamondsuit =  \spadesuit - \frac{1}{2}\id_{A^*}\clubsuit=0 $, i.e. if and only if $ \spadesuit = \clubsuit = \diamondsuit=0 $, as was to be shown.
\end{pf}

We call the pair  $(\mu^0_{\phi} , \nu^{1}_{X}) $
the \emph{Poissonified Lie bialgebroid} of the Jacobi bialgebroid $(\mu+\phi,\nu+X) $.

\begin{rmk}
The last corollary allows one to resume all the conditions listed in theorem \ref{theo:main} to the single condition $\bb{ \mu^0_{\phi} + \nu^{1}_{X},\mu^0_{\phi} + \nu^{1}_{X} }=0$,
or, more generally $\bb{ \mu^c_{\phi} + \nu^{1-c}_{X},\mu^c_{\phi} + \nu^{1-c}_{X} }=0$ for some $c \in \R$.
Said differently, a pair $(\mu+\phi,\nu+X)$ of pre-Jacobi algebroids is a Jacobi bialgebroid
if and only if $ \bb{ \mu^0_{\phi} + \nu^{1}_{X},\mu^0_{\phi} + \nu^{1}_{X} }=0$.
\end{rmk}

\section{Poisson-Jacobi manifolds and its Jacobi bialgebroids. }
\label{sec:PoissonJacobiMfd}

\noindent
{\bf Notations}
Throughout this section, we shall use the following conventions.

Let $B \to M$ be a vector bundle endowed with a Lie algebroid $\mu_B \in \superF{B}$. As in section \ref{sec:prelim}-B, we call $A=B \oplus \R \to M $ the direct sum of $B$  with the trivial bundle $\R \times M \to M$,  $\epsilon$  the section of $A:= B \oplus \R \to M$ given by $m \to (0_m,1) \in A_m = B_m \oplus \R$, and $\phi$ the section of  $A^*$ given by $\phi (b+f\epsilon) := f$ for every $b \in \Gamma(B), f \in \FF(M)$. As mentioned in section \ref{sec:prelim}-B, there is a natural inclusion $\superF{B} \subset \superF{A}$: since it preserves the big bracket, $\mu_B$ can also be considered as a Lie algebroid on $A \to M$. To avoid confusion, we shall denote by $\mu$ this Lie algebroid on $A$.

 It is clear that $\phi \in \Gamma(A^*)$ is a cocycle for the Lie algebroid structure $\mu$,
so that $ \mu+\phi$ is a Jacobi algebroid on $A \to M$. By lemma \ref{lem:alg_associadas_Jaco_alg} therefore, $ \mu_\phi^0= \mu + \phi \, \partialt  $ is a Lie algebroid on $p^* A \to (M \times \R)$,
where $\partialt \in \superF{p^*A}$ is as in section \ref{sec:prelim}-B, and $t\in \FF(M \times \R) $ stands for the parameter on $\R $ (as in section \ref{sec:prelim}-A).

\begin{rmk}
Explicitly, the Lie algebroid associated to $\mu$  on $A= B \oplus \R$ is the Lie algebroid
(which appears in \cite{IM}) with bracket and anchors:
  \begin{eqnarray*} \lb{b_1+f_1 \epsilon ,b_2+f_2 \, \epsilon}_A &:=& [b_1,b_2]_B+\rho_A(b_1)[f_2] \, \epsilon- \rho_A(b_2)[f_1]\, \epsilon  \\
   \hbox{ and } \rho_A ( b,f) &:=& \rho_B (b) \quad \quad \forall b_1,b_2,b \in \Gamma(B), f_1,f_2,f \in \FF(M) .
   \end{eqnarray*}
 \end{rmk}

\subsection{Poissonization}

By a \emph{Jacobi structure on the Lie algebroid $\mu_B$ on $B \to M$}, we mean a pair $(\pi,E)$ with $\pi \in \Gamma(\wedge^2 B)$
and $E \in \Gamma(B) $ such that:
     $$ \lb{\pi,\pi}_{\mu_B}= -2 E \pi \hbox{ and } \lb{E,\pi}_{\mu_B} =0.$$
Of course, the coefficient $-2$ is arbitrary, and could be turned into $1$ by replacing $E$ by~$-2 E$.

\begin{lem}\label{lem:Peixe}
We use the notations introduced in the beginning of section \ref{sec:PoissonJacobiMfd}.
For every $\pi \in \Gamma(\wedge^2 B)$ and $E \in \Gamma(B)$, the following are equivalent:
\begin{enumerate}
\item[(i)] $(\pi,E)$ is a Jacobi structure for the Lie algebroid $\mu_B$ on $B$,
\item[(ii)] $ P_{\pi,E}:= e^{-t} (\pi + \epsilon E) $ is a Poisson structure for the Lie algebroid $\mu_\phi^0$ on $ p^* A \to (M \times \R) $.
\end{enumerate}
\end{lem}
\begin{pf}
Recall that, by definition, $P_{\pi,E} \in \Gamma( \wedge^2 p^* A)$ is a Poisson structure for the Lie algebroid $\mu_\phi^0 $ if and only if $ \lb{P_{\pi,E},P_{\pi,E}}_{\mu_\phi^0 } =0 .$
A direct computation gives the following:
 \begin{eqnarray*} \lb{P_{\pi,E},P_{\pi,E}}_{\mu_\phi^0 } &=&\bb{ \bb{P_{\pi,E}, \mu_\phi^0 } , P_{\pi,E}} \\
  &=& \bb{ \bb{ e^{-t}(\pi + \epsilon E) , \mu + \phi \partialt } ,e^{-t}(\pi + \epsilon E)}  \\
   & =&  \bb{  e^{-t} \bb{ \pi,\mu} + e^{-t} \epsilon  \bb{ E,\mu} + \bb{e^{-t}, \partialt } \phi (\pi + \epsilon E) -e^{-t} E \partialt,e^{-t}(\pi + \epsilon E)} \\
   & =&  \bb{  e^{-t} \bb{ \pi,\mu} + e^{-t} \epsilon  \bb{ E,\mu}  - e^{-t} \phi (\pi + \epsilon E) -e^{-t} E\partialt, e^{-t}(\pi + \epsilon E)} \\
  &=& e^{-2t} \bb{  \bb{ \pi,\mu} +  \epsilon  \bb{ E,\mu}  + \phi (\pi + \epsilon E) - E \partialt, (\pi + \epsilon E)} \\
  & & \quad +e^{-t}\bb{ -E \partialt , e^{-t}}(\pi+\epsilon E)\\
  &=& e^{-2t} (\lb{\pi,\pi}_{\mu_B} + 2 E \pi+  2\epsilon \lb{E,\pi}_{\mu_B} ) .\\
   \end{eqnarray*}
Now, a section of $\wedge^2 A$  of the form $\Theta_1 + \epsilon \Theta_2 =0$, with $\Theta_1 \in \Gamma(\wedge^3 B)$ and
$\Theta_2 \in \Gamma(\wedge^2 B)$ vanishes if and only $\Theta_1=\Theta_2=0$, so that $  \lb{P_{\pi,E},P_{\pi,E}}_{\mu_\phi^0 } =0 $ if and only if
$ \lb{\pi,\pi}_{\mu_B} + 2E \pi= \lb{E,\pi}_{\mu_B} =0 $.
\end{pf}

\begin{numrmk}\label{rmk:poissonization}
In general \cite{GM,IM}, the Poissonization is constructed on the direct product of the Lie algebroid $B \to M $ with the tangent algebroid $T\R \to \R$. But this Lie algebroid is indeed canonically isomorphic
to $p^* (B \oplus \R) \to (M \times R)$, and our construction matches the usual one.

Moreover, when $B=TM$ is the Lie algebroid of vector fields, $p^* A$ is isomorphic to $T(M \times \R)$,
and, under the previous isomorphism, $\epsilon $ is the vector field $\frac{\partial }{\partial t}$.
Specializing lemma \ref{lem:Peixe}, one obtains that a pair $(\pi,E)$
 (with $\pi$ a bivector field and $E$ a vector field) is a Jacobi structure on $TM$ if and only if $
e^{-t}(\pi + \frac{\partial }{\partial t} \wedge E )$ is a Poisson structure on $M \times \R$.
\end{numrmk}

\subsection{The Jacobi structure of a Jacobi bialgebroid}

Exactly like Lie bialgebroids induce Poisson structures on their underlying manifold,
Jacobi bialgebroids induce Jacobi structures on their underlying manifolds, see \cite{IM}.
Our formalism gives quite easy descriptions of these structures.

\begin{prop}
Let $(\mu+\phi,\nu+X)$ be a Jacobi bialgebroid on a vector bundle $A \to M$,
then
\begin{enumerate}
\item The assignment
$$\begin{array}{rcl}
 \FF(M)&\to &\FF(M) \\
 f &\mapsto &\bb{\bb{\mu,X},f}
\end{array} $$
is a vector field $\underline{E}$ on $M$. This vector fields is also given by the assignment
$f \mapsto-\bb{\bb{\nu,\phi},f} $.
\item The assignment
$$\begin{array}{rcl}
\FF(M) \times \FF(M)&\to &\FF(M) \\
 (f,g) &\mapsto &\bb{\bb{\mu,f},\bb{\nu,g}}
\end{array}
$$
is a bivector field $\underline{\pi}$ on $M$.
\item The pair $(\underline{\pi},-\underline{E})$ is a Jacobi structure on $TM$.
\end{enumerate}
\end{prop}
\begin{pf}
It is clear that both assignments are derivations in each variables for degree reasons.
Let us check also the second one is skew-symmetric. A direct computation gives, for every $f,g \in \FF(M)$, that $\underline{\pi}[f,g]-\underline{\pi}[g,f]= \bb{ \bb{ \bb{ \mu,\nu }, f} , g } $.
According to theorem \ref{theo:main}, the quantity $\spadesuit=0$, so that
 the previous can be rewritten as:
$$ \underline{\pi}[f,g]-\underline{\pi}[g,f] = -\bb{ \bb{\mu X+\nu\phi - \id_{A}\left(  -\phi X + \frac{\bb{ \mu,X }- \bb{ \nu,\phi }}{2}\right),\hbox{$f$}},g  }   ,$$
a quantity that vanishes for degree reasons. This proves the first two items.

We now turn our attention to the third item.
According to theorem \ref{cor:ligacao_com_bialg}, the pair $(\mu^0_\phi, \nu^1_X)$ is a Lie bialgebroid. Hence, according to proposition 3.6 in \cite{MX}, the assignment $\FF (M \times \R) \times \FF(M \times \R)
\mapsto \FF(M \times \R) $ given by
   $$ (f,g) \mapsto \bb{ \bb{ \mu^0_\phi , f } , \bb{ \nu^1_X , g }}  $$
is a Poisson structure on $M \times \R$, that we denote by $P$. In view of lemma \ref{lem:Peixe} (more precisely remark \ref{rmk:poissonization}), in order to complete the proof it suffices to check that $P$ is precisely the Poissonization
of the pair $(\underline{\pi},-\underline{E})$. The computation goes as follows for every two functions $f,g$
on $M\times \R $:
\begin{eqnarray*}
 P[f,g] &=& \bb{\bb{ \mu^0_\phi , f},\bb{ \nu^1_X, g} } \\
 &=& \bb{\bb{ \mu + \partialt \phi , f},\bb{e^{-t}(\nu + (\partialt + \id_{A^*} )X) , \hbox{$g$}}}  \\
 &=& \bb{\bb{ \mu ,f} + \phi \frac{\partial f}{\partial t},e^{-t}\left(\bb{\nu,g} + X \frac{\partial g}{\partial t}\right)   }  \\
 &=& e^{-t} \left( \bb{\bb{ \mu ,f},\bb{\nu,g}  } + \bb{X,\phi}  \frac{\partial f}{\partial t}\frac{\partial g}{\partial t}
 +   \bb{\bb{ \mu ,f} ,   X  \frac{\partial g}{\partial t}}+
 \bb{\phi \frac{\partial f}{\partial t}, \bb{\nu,g}}  \right) \\
 &=& e^{-t} \left( \bb{\bb{ \mu ,f},\bb{\nu,g}  }
 +   \bb{\bb{ \mu ,f} ,   X \frac{\partial g}{\partial t}}+
 \bb{\phi \frac{\partial f}{\partial t}, \bb{\nu,g}}  \right) \\
 & & \quad \quad \quad \hbox{ (since $\bb{\phi,X}=\diamondsuit=0$ )} \\
&=& e^{-t} \big(\bb{\bb{ \mu ,f},\bb{\nu,g}  }
 +   \bb{\bb{ \mu ,f} ,   X} \frac{\partial g}{\partial t}+
  \frac{\partial f}{\partial t}\bb{\phi, \bb{\nu,g}}  \big) \\
&=& e^{-t} \big(\underline{\pi}[f,g]
 +  \underline{E}(f) \frac{\partial g}{\partial t}-
 \frac{\partial f}{\partial t}\underline{E}(g)  \big)    \\
&=&  e^{-t} ( \underline{\pi} + \frac{\partial}{\partial t} \wedge (-\underline{E})  ) [f,g].
\end{eqnarray*}
This completes the proof.
\end{pf}

\subsection{The Jacobi bialgebroid of a Jacobi structure}

For every Poisson structure $P$ on a Lie algebroid $\nu$ defined on a vector bundle $C$, the function $\bb{P,\nu}$
is a Lie algebroid structure on $C^*$. Moreover, the pair $(\nu,\bb{P,\nu})$ is a Lie bialgebroid. We refer to \cite{MX} for these classical results .
We call the Lie algebroid  $\bb{P,\nu}$ (resp. the Lie bialgebroid $(\nu,\bb{P,\nu})$) the \emph{Lie algebroid (resp. the Lie bialgebroid) associated to the Poisson bivector $P$}.
In particular, it follows from lemma \ref{lem:Peixe} that, for every Jacobi structure $(\pi,E) $ on a Lie algebroid $\mu_B$ on $B \to M $, the function in $\superF{p^* A}$ (see section \ref{sec:prelim} for the notations)
 defined by
  $$\nu_{\mu,\pi,E} := \bb{  P_{\pi,E} , \mu^0_\phi } = \bb{  P_{\pi,E} , \mu + \phi \partialt } $$
  is a Lie algebroid structure on $(p^* A)^* \to (M \times \R) \simeq p^* A^* \to (M \times \R)$.
  Moreover, the pair $(\mu^0_\phi,\nu_{\mu,\pi,E} )$ is a Lie bialgebroid on the vector bundle
  $p^* A \to (M \times \R) $.

We start with a lemma:

\begin{lem}\label{lem:poissoniacao}
We use the notations introduced in the beginning of section \ref{sec:PoissonJacobiMfd}.
Let $(\pi,E)$ be a Jacobi structure on the algebroid $\mu_B$ on $B\to M$.
Define
 $$\nu:=\bb{ \pi, \mu}+ \phi \pi+ \epsilon \bb{\mu,E} -E \, \id_B  \in \superF{A}^{1,2}.$$
The Poissonization of the pre-Jacobi algebroid $\nu-E$ on $A$ is the Lie algebroid
$\nu_{\mu,\pi,E}:= \bb{ P_{\pi,E},\mu^0_\phi }  $ associated with the Poisson structure $P_{\pi,E} $.
\end{lem}
\begin{pf}
We have to show that the functions ${\nu}_{-E}^{1} = e^{-t}(\nu- E(\partialt + \id_{A^*} ))  $
and $ \nu_{\mu,\pi,E} = \bb{ P_{\pi,E} , \mu + \phi \partialt }  $ are equal.
This follows from a comparison between
\begin{eqnarray*}
\nu_{\mu,\pi,E} & = & \bb{ e^{-t} (\pi + \epsilon E) ,\mu + \phi \partialt } \\
 & =& e^{-t} (\bb{   (\pi + \epsilon E),\mu} + \phi (\pi + \epsilon E ) - \partialt E ) \\
  & =& e^{-t} (\bb{\pi,\mu} + \bb{\epsilon E,\mu} + \phi (\pi + \epsilon E ) - \partialt E ) \\
   & =& e^{-t} (\bb{\pi,\mu} + \epsilon \bb{E,\mu} + \phi (\pi + \epsilon E ) - \partialt E )
\end{eqnarray*}
and
\begin{eqnarray*}
\nu^1_{-E}
 &= &e^{-t}(\nu+E(\partialt - \id_A )) \\
 &= & e^{-t}(\bb{ \pi, \mu}+ \phi \pi+ \epsilon \bb{\mu,E} -E \, \id_B - E(\partialt - \id_A )) \\
 & =& e^{-t}(\bb{\pi, \mu}+ \phi \pi +\epsilon \bb{\mu,E}-  \partialt E + (\id_A -id_B ) E)\\
 & =& e^{-t}( \bb{\pi, \mu}+ \phi \pi + \epsilon \bb{E,\mu}- \partialt E + \phi \epsilon E)
\end{eqnarray*}
where, in the last line, we have used the relation $\id_A= \id_B + \phi \epsilon$.
\end{pf}

\begin{rmk}
A direct comparison shows that the bracket and anchor of the associated pre-Lie algebroid $\nu$ on $A \to M$ are precisely the bracket and anchor of the Lie algebroids that appears in \cite{Ker,GM,IM}.
\end{rmk}

The main result of this section is an immediate consequence of lemma
\ref{lem:alg_associadas_Jaco_alg} and of lemma~\ref{lem:poissoniacao}.

\begin{prop}
We use the notations introduced in the beginning of section \ref{sec:PoissonJacobiMfd}.
Let $(\pi,E)$ be a Jacobi structure on the algebroid $\mu_B$ on $B\to M$. Then,
 \begin{enumerate}
 \item the pair $( \mu+\phi,\nu - E) $ is a Jacobi bialgebroid structure,
 where:
 $$\nu:=\bb{ \pi, \mu}+ \phi \pi+ \epsilon \bb{\mu,E} -E \, \id_B  \in \superF{A}^{2,1}.$$
\item The Poissonization of this  structure is the Lie bialgebroid
$(\mu^0_\phi,\nu_{\mu,\pi,E} ) $ associated to the Poisson structure $P_{\pi,E}$
(i.e. the Poissonization of the Jacobi structure $(\pi,E)$).
\end{enumerate}
\end{prop}
We recapitulate the results of this section in a commutative diagram.
Let $\mu_B$ be a Lie algebroid on a vector bundle $B$, equipped with a Jacobi structure $(\pi,E)$.
Recall that $A$ stands for $B \oplus \R$ (and $\mu$ is the natural extension of $\mu_B$ to that bundle), while $p^*A $ stands for the pull-back of $A \to M$ on $M \times \R$ (and can be equipped with the Lie algebroid structure $\mu_\phi^0:=\mu+\partialt \phi$).
Let $P_{\pi,E}$  be the Poissonization of the Jacobi structure structure $(\pi,E)$ (defined on the Lie algebroid $\mu_\phi^0$
on $p^* A$). Last, let $\nu := \bb{ \pi, \mu}+ \phi \pi+ \epsilon \bb{\mu,E} -E \, \id_B $ be the Jacobi algebroid structure on $A$ defined above. The following diagram commutes
$$\xymatrix{
*\txt{$(\pi, E)$\\{\hspace*{1em}Jacobi structure \hspace*{1em}}\\ on $B$ for $\mu_B$}\ar[ddd] \ar[rr]^{\textrm{Poissonization}} &\quad\quad\quad\quad\quad &*\txt{ $P_{\pi,E}$\\{\hspace*{1em}Poisson bivector\hspace*{1em}}\\ on $p^*A$ for $\mu^0_{\phi}$}\ar[ddd]\\
&\quad\quad\quad\quad\quad&\\
&\quad\quad\quad\quad\quad&\\
*\txt{$(\mu+\phi, \nu - E)$\\ {\hspace*{1em}Jacobi bialgebroid\hspace*{1em}}\\ on $(A, A^*)$} \ar[rr]^{\textrm{Poissonization}} &\quad\quad\quad\quad\quad & *\txt{ $\left(\mu^0_{\phi}, \nu^1_{-E}\right)$\\{\hspace*{1em}Lie bialgebroid\hspace*{1em}}\\ on $(p^*A, p^*(A^*))$}}$$

\section{Quasi-Jacobi bialgebroids}

In \cite{AJR}, quasi-Jacobi bialgebroids are introduced, as follows:

\begin{defi}
Let $A \to M$ be a vector bundle. A quasi-Jacobi bialgebroid is a
pre-Jacobi bialgebroid $(\mu+\phi,\nu+X) $, together with a section $\ZZZ$ of $\wedge^3 A $,
such that $ \mu+\phi$ is indeed a Jacobi algebroid, and such that the following compatibility conditions are satisfied
    \begin{equation}\label{cond_Joana}
    \diff_{\nu,X} (\ZZZ)     \hbox{ and }   (\diff_{\nu,X})^2 (P) =\lb{-\ZZZ,P}_{\mu,\phi} =  0 \quad \hbox{ for all $P \in \Gamma(\bigwedge A)$}
    \end{equation}
where $\LB_{\mu,\phi},\diff_{\nu,X}$ are as in definition \ref{def:jacobi_algebroids}.
\end{defi}

As a slight modification of theorem \ref{theo:main}, we obtain:

\begin{theo}\label{theo:quasi_main}
Let $A \to M $ be a vector bundle, $\mu+\phi$ a pre-Jacobi algebroid structure on $A$, $\nu+X$ a pre-Jacobi algebroid structure on $A^*$, and $\ZZZ$ a section of $\wedge^3 A$. Define
$\diamondsuit $, $\clubsuit $ and $ \spadesuit$ as in (\ref{eq:horrible}).
The triple $(\mu+\phi,\nu+X,\ZZZ) $ is a quasi-Jacobi bialgebroid if and only if
 $$ \diamondsuit=\clubsuit=\spadesuit= \bb{\mu+\phi,\mu+\phi} =\bb{\nu,\ZZZ }+X \ZZZ = \bb{\nu,X} + \bb{\ZZZ,\phi} = \maltese =0 ,$$
 with $ \maltese = \frac{1}{2} \bb{\nu,\nu}+ \bb{\ZZZ,\mu} + \id_A \bb{\ZZZ,\phi} + 2 \ZZZ \phi $.
\end{theo}
\begin{pf}
According to theorem \ref{theo:main}, the pair $(\mu+\phi,\nu+X) $ is a pre-Jacobi bialgebroid
if and only if $ \diamondsuit=\clubsuit=\spadesuit= 0$. Also, $ \bb{\mu+\phi,\mu+\phi}=0$ if and only if $\mu+\phi$ is a Jacobi algebroid. The condition $\diff_{\nu,X} (\ZZZ) =0$
is tantamount to $\bb{\nu,\ZZZ }+X \ZZZ=0$ in view of the expression of the differential given in definition \ref{def:jacobi_algebroids}.

Spelling out the condition $(\diff_{\nu,X})^2 (P) =\lb{-\ZZZ,P}_{\mu,\phi}$ with the help of the differentials and brackets introduced in definition \ref{def:jacobi_algebroids} gives:
  $$ \bb{ \nu, \bb{\nu,P} +XP  } + X(  \bb{\nu,P} +XP   )= -\bb{\bb{Z ,\mu+\id_A \phi } ,P }
  + \bb{ Z \phi ,P }-\bb{Z,\phi} P, $$
  which can be rewritten as
   $$ \bb{ \maltese , P} + ( \bb{\nu,X} +\bb{Z,\phi})P=0 $$
This condition is therefore satisfied if $\bb{\nu,X} + \bb{\ZZZ,\phi} = \maltese =0$.
Conversely, if this condition is satisfied for $P=1$, then $\bb{\nu,X} + \bb{\ZZZ,\phi} =0$. If the condition  is satisfied for all $P\in \Gamma(\bigwedge A)$ therefore, we have $ \bb{\maltese,P}=0$ for all $P\in \Gamma( A)$, which, in view of the seventh item listed in section \ref{sec:defi}, since $\maltese \in \superF{A}^{3,1}$, implies that $\maltese =0$. This completes the proof.
\end{pf}

Recall \cite{Roytenberg} that a quasi-Lie bialgebroid on a vector bundle $B$ is a
triple of  functions in, respectively, $  \superF{A}^{1,2}, \superF{A}^{2,1} , \superF{A}^{3,0} $, and whose sum
commutes with itself. Theorem \ref{theo:quasi_main} admits the following corollary, to be compared with corollary \ref{cor:ligacao_com_bialg}:

\begin{coro}\label{cor:ligacao_com_quasi_bialg}
Let $A \to M $ be a vector bundle, $\mu+\phi$ a pre-Jacobi algebroid structure on $A$, $\nu+X$ a pre-Jacobi algebroid structure on $A^*$,
and $\ZZZ \in \Gamma(\wedge^3 A)$.
 The following are equivalent:
\begin{enumerate}
\item[(i)] the triple $(\mu+\phi,\nu+X,\ZZZ  ) $ is a quasi-Jacobi bialgebroid on $A$;
\item[(ii)] the triple $(\mu^0_{\phi} , \nu^{1}_{X}, e^{-2t} \ZZZ) $ is a quasi-Lie bialgebroid on $p^* A$.
\end{enumerate}
(In the previous, $ \mu^0_{\phi} , \nu^{1} _{X}$ are the pre-Lie algebroids on $p^*A \to (M \times \R)$ constructed in lemma \ref{lem:alg_associadas_Jaco_alg} and $t \in \FF(M \times \R)$ is, as in section \ref{sec:prelim}-A,  the projection on the second component)
\end{coro}
\begin{pf}
The triple $ (\mu^0_{\phi} , \nu^{1} _{X}, e^{-2t} \ZZZ) $ is a quasi-Lie bialgebroid on $p^* A$
if and only if $$\bb{\mu^0_{\phi} + \nu^{1} _{X}+ e^{-2t} \ZZZ, \mu^0_{\phi} + \nu^{1} _{X}+ e^{-2t} \ZZZ}=0.$$
For degree reasons, this condition splits into the four conditions:
 $$
 \begin{array}{rclr}
 \bb{\mu^0_{\phi},\mu^0_{\phi}}  & = & 0 & (\hbox{*})\\
 \bb{\mu^0_{\phi},\nu^1_{X}}    & = & 0  & (\hbox{**})\\
 \bb{\nu^1_{X},\nu^1_{X}}   & = & -2 \bb{ \mu^0_{\phi} , e^{-2t} \ZZZ} & (\hbox{***})\\
\bb{\nu^1_{X},e^{-2t} \ZZZ}   & = & 0 & (\hbox{****})\\
 \end{array}
 $$
As follows from lemma \ref{lem:alg_associadas_Jaco_alg}, condition (*) holds if and only if $ \bb{\mu+\phi,\mu+\phi}=0 $.
As follows from item \emph{(1)} in corollary \ref{cor:ligacao_com_bialg},
condition (**) holds if and only if $(\mu+\phi,\nu+X)$ is a pre-Jacobi algebroid, i.e.,
in view of theorem \ref{theo:main}, if and only if $ \diamondsuit=\clubsuit=\spadesuit=0$.
As follows from a direct computation, condition (****) holds if and only if
 $$ \bb{\nu^1_{X},e^{-2t} \ZZZ} = e^{-3t}( \bb{\nu,\ZZZ }+2\ZZZ X -3 \ZZZ X ) =0, $$
 i.e. if and only if $ \bb{\nu,\ZZZ }+X \ZZZ =0$.

In view of theorem \ref{theo:quasi_main} therefore, we are therefore left with the task of showing that  condition (***) holds if and only if $ \bb{\nu,X} + \bb{\ZZZ,\phi} = \maltese=0$.
We check it as follows. First, we compute:
 \begin{eqnarray*}
 \bb{\nu^1_{X},\nu^1_{X}}&=& \bb{ e^{-t}(\nu+ X (\partialt + \id_{A^*}  )) ,e^{-t}(\nu+ X  (\partialt + \id_{A^*} ))} \\
 &=&   e^{-2t}( \bb{\nu,\nu} +  2\bb{\nu,X} \id_{A^*}   +2 \bb{\nu,X} \partialt) .
\end{eqnarray*}
Second, a direct computation gives $\bb{ \mu^0_{\phi} , e^{-2t} \ZZZ}=e^{-2t} (  \bb{\mu,\ZZZ}-2 \phi \ZZZ + \partialt \bb{\phi,Z}) $, so that, reordering terms, we obtain that (***) is equivalent to
$$ \maltese + \partialt (\bb{\nu,X}+ \bb{Z,\phi} ) = 0 , $$
which is itself equivalent to the vanishing of both $ \maltese$ and $ \bb{\nu,X}+ \bb{Z,\phi}$ (as noticed in section \ref{sec:prelim}-A). The result follows.
 \end{pf}

\end{document}